\documentclass[reqno,11pt]{amsart}
\usepackage{amssymb, amsthm}
\usepackage[numbers]{natbib}
\usepackage[all]{xy}
\input epsf
\usepackage{lscape}

\numberwithin{figure}{section}

\theoremstyle{plain}
\newtheorem{thm}{Theorem}[section]
\newtheorem{lem}[thm]{Lemma}

\newtheorem{conj}[thm]{Conjecture}

\theoremstyle{definition}
\newtheorem{exl}[thm]{Example}

\title{A note on the topological sliceness of some 2-bridge knots}

\author{Allison N. Miller} \thanks{The author was supported by NSF grant DMS-
1148490.}
\address{Department of Mathematics, University of Texas, Austin, TX 78712, USA}

\begin{document}

\begin{abstract}
We use twisted Alexander polynomials to show that certain algebraically slice 2-bridge knots are not topologically slice, even though all prime power Casson-Gordon signatures  vanish. We also provide some computations indicating the efficacy of Casson-Gordon signatures in obstructing the smooth sliceness of 2-bridge knots. 
\end{abstract}
\bibliographystyle{alpha}

\maketitle

\section{Introduction}

Although 2-bridge knots $K_{p,q}$ are generally well understood, their algebraic and topological slice status is not. One of the only easily applicable statements in terms of $p$ and $q$ is that if $K_{p,q}$ is algebraically slice then $|H_1(\Sigma_{2}(K_{p,q}))|=p$ must be a square.
In \cite{CG86}, Casson and Gordon gave the first examples of algebraically slice knots which were not ribbon, smoothly slice, or even topologically slice. For an algebraically slice knot $K$, every prime-power branched cover $\Sigma_{p^n}(K)$ has first homology with order equal to some square $m^2$. For any $k$ dividing $m$ and $r$ with $0 \leq r \leq k-1$, there is a Casson-Gordon signature $\sigma_{CG}(K; p^n, k,r)$. 
If $K$ is ribbon, then $\sigma_{CG}(K; p^n, k,r)$ must vanish for all choices of $p^n, k,$ and $r$ as above; however, sliceness (smooth or topological)  only implies that these signatures must vanish for $k$ a  prime power. 
The signatures associated to the double branched cover of rational knot $K_{m^2,q}$ are particularly computable; in fact, there is a combinatorial formula in terms of counts of integer points in triangles. Casson and Gordon observed in \cite{CG86} that the only known rational knots for which all $\sigma_{CG}(K; 2, k, r)$ vanished belonged to a certain family $\mathcal{R}$ of ribbon knots. 

\begin{conj}[\cite{CG86}, \cite{EL09}]\label{cgsig}
Suppose $K_{m^2,q}$ is a 2-bridge knot. Then $K_{m^2,q}$ is ribbon if and only if 
 all Casson-Gordon signature invariants associated to the double branched cover vanish if and only if $K_{m^2, q}$ is in $\mathcal{R}$. \end{conj}

Lisca partially resolved this question by classifying the smooth sliceness of rational knots. 
\begin{thm}[\cite{Lisca07}]
$K_{p,q}$ is smoothly slice if and only if $K_{p,q}$ is ribbon if and only if $K_{p,q} \in \mathcal{R}.$
\end{thm}
Despite this classification, the question of exactly when the Casson-Gordon signature invariants vanish remains open.\footnote{See \cite{EL09} for more discussion of 
Conjecture \ref{cgsig} from a number-theoretic perspective.} Answering this question would give additional information about which 2-bridge knots are topologically slice.  In particular, an affirmative answer would show that for $m$ is a prime power the topological sliceness, smooth sliceness, and ribbonness  of $K_{m^2,q}$ all coincide with the vanishing of the double branched cover Casson-Gordon signature invariants.

The first algebraically slice 2-bridge knot for which the Casson-Gordon signature invariants do not obstruct sliceness is $K_{225, 94}$, as observed in \cite{CG86}. We compute a twisted Alexander polynomial associated to the double branched cover and observe that this polynomial shows that  $K$ is not topologically slice. We also give some computations indicating the effectiveness of the Casson-Gordon signature invariants (particularly when combined with the classical Alexander polynomial) at obstructing the topological sliceness of $K_{m^2,q}$ for small values of $m$.

\section{Twisted Alexander polynomials}
In general, twisted homology and twisted Alexander polynomials can be defined for spaces $Y$ which are homotopy equivalent to finite CW complexes.\footnote{We follow the much more thorough exposition of \cite{KL99a} and \cite{HKL10}.}
Let $\tilde{Y}$ denote the universal cover of $Y$, so $C_*(\tilde{Y})$ is acted on by the left by $\pi=\pi_1(Y)$. 
Given $M$ a $(\mathbb{F}[t^{\pm1}], \mathbb{Z}[\pi])$ bimodule, the $M$-twisted chain complex of $Y$ is  $C_*(Y, M):= M \otimes_{\mathbb{Z}[\pi]} C_*(\tilde{Y})$. 
Note that $C_*(Y,M)$ and hence $H_k(Y,M)= H_k(C_*(Y,M))$ inherit a left $\mathbb{F}[t^{\pm1}]$-module structure from $M$. The \emph{twisted Alexander polynomial} $\Delta_{Y,M}(t)$ associated to $Y$ and $M$ is defined to be the order of $H_1(Y,M)$ as a $\mathbb{F}[t^{\pm1}]$ module. 

Let $K$ be a knot,  $X$ denote its exterior, $X_m$ denote the canonical cyclic $m$-fold cover of $X$, and $\Sigma_m$ denote the corresponding branched cover of $S^3$ over $K$.  
There is a canonical map $\epsilon: \pi_1(X) \to \mathbb{Z}$. Let  $\epsilon_m$ be the composition $\pi_1(X_m) \hookrightarrow \pi_1(X) \xrightarrow{\epsilon} \mathbb{Z}$ restricted to its image. 
Choose $n$ a prime power dividing $|H_1(\Sigma_m)|$, a map $\chi:H_1(X_m) \to H_1(\Sigma_m) \to \mathbb{Z}_n$, and $\xi_n$ a primitive $n^{th}$ root of unity. 
Then $M= \mathbb{Q}(\xi_n)[t^{\pm1}]$  has a $(\mathbb{Q}(\xi_n)[t^{\pm 1}], \mathbb{Z}[\pi_1(X_m)])$-bimodule structure given by polynomial multiplication on the left and  $\mathbb{Z}[\pi_1(X_m)]$ action defined by 
$ p(t) \cdot \gamma= \xi_n^{\chi(\gamma)} t^{\epsilon_m(\gamma)} p(t)$ for $\gamma \in \pi_1(X_m)$.\footnote{Note that we often abuse notation by blurring the distinction between an element of a fundamental group and its image in first homology.}
It is often convenient to consider the \emph{reduced} twisted Alexander polynomial
 $ \widetilde{\Delta}_{X, M}(t):= \Delta_{X, M}(t) (t-1)^{-s}, \text{ where } s=0 \text{ if } \chi$ is trivial and $s=1 \text{ else. }$
These metabelian twisted Alexander polynomials $\Delta_{X_m, M}$ give an obstruction to the topological sliceness of $K$, as follows. 

\begin{thm}[\cite{KL99a}]\label{notslice}
Let $K$ be a topologically slice knot and $a,b$ distinct primes with $b \neq 2$.
Let  $m=a^r, n=b^s$. 
 Then there exists an invariant metabolizer $N \leq H_1(\Sigma_m)$ such that if $\chi: H_1(X_m) \to H_1(\Sigma_m) \to \mathbb{Z}_n$  vanishes on $N$ then the corresponding reduced twisted Alexander polynomial  is a norm in $\mathbb{Q}(\xi_n)[t^{\pm1}]$.  
 That is, there exists $\lambda \in \mathbb{Q}(\xi_n)$, $k \in \mathbb{Z}$, and $f(t) \in \mathbb{Q}(\xi_n)[t^{\pm1}]$ such that $\widetilde{\Delta}_{X_m, \mathbb{Q}(\xi_n)[t^{\pm1}]}(t)= \lambda t^k f(t) \overline{f(t^{-1})}$. 
\end{thm}

Note that when $K=K_{p,q}$ is 2-bridge and $m=2$ the application of Theorem \ref{notslice} is particularly straightforward, since $\Sigma_2(K_{p,q})=L_{p,q}$. 
Let $k$ be a prime dividing $p$. As a $\mathbb{F}_k[\mathbb{Z}_2]$ module, $H_1(\Sigma_2, \mathbb{Z}_k)$ must be isomorphic to the direct sum of modules of the form $\mathbb{F}_k[t]/f(t)$ , where $f(t)$ divides both $\Delta_K(t)$ and $t^2-1$ in $\mathbb{F}_k[t]$.
So $H_1(\Sigma_2, \mathbb{Z}_k) \cong \left(\mathbb{F}_k[t]/\langle t+1 \rangle \right)^r$. However, since $\Sigma_2$ is a lens space, 
the first homology $H_1(\Sigma_2) \cong \mathbb{Z}_p$ is cyclic. 
So $r=1$ and $H_1(\Sigma_2, \mathbb{Z}_k) \cong \mathbb{F}_k[t]/\langle t+1 \rangle$ is an irreducible $\mathbb{F}_k[\mathbb{Z}_2]$ module. 
Therefore, as observed by \cite{HKL10}, any metabolizer $N \leq H_1(\Sigma_2)$ must have trivial image $ \overline{N} \leq H_1(\Sigma_2, \mathbb{Z}_k)$. 
In order to obstruct the topological sliceness of $K_{p,q}$, it therefore suffices to show that a single twisted Alexander polynomial coming from a character factoring through $H_1(\Sigma_2(K_{p,q}), \mathbb{Z}_k)$ is not a norm.

Computation of the twisted Alexander polynomials of covers is significantly simplified by Herald, Kirk, and Livingston's reinterpretation in terms of certain twisted Alexander polynomials corresponding to more complicated representations of the base space. In this context, their work in \cite{HKL10} gives the following. 
Let $H=H_1(\Sigma_2, \mathbb{Z}_k)= \mathbb{F}_k[t]/\langle t+1 \rangle$, so $\mathbb{Z} \ltimes H$  has multiplication given by $(x^i,v) \cdot (x^j, w)= (x^{i+j}, t^{-j}\cdot v+w)= (x^{i+j}, (-1)^{-j}v +w).$ 
Choose a meridian $\mu \in \pi_1(X)$ with $\epsilon(\mu)=1$. 
Then there is a correspondence between equivariant\footnote{
Note that conjugation by $\mu$ gives an automorphism of $\pi_1(X_2) \leq \pi_1(X)$, and $\rho$ is \emph{equivariant} if
$\rho(\mu \gamma \mu^{-1}) = t \cdot \rho(\gamma)$ for any $\gamma \in \pi_1(X_2)$ and $\mu$ our preferred meridian.}
 homomorphisms $\rho: \pi_1(X_2) \to H$ and homomorphisms $\tilde{\rho}: \pi_1(X) \to \mathbb{Z} \ltimes H$ that extend  $\epsilon|_{\pi_1(X_2)} \times \rho: \pi_1(X_2) \to 2\mathbb{Z} \times H$ and with $\tilde{\rho}(\mu)=(x,0)$.\footnote{
 Given $\rho$, this correspondence associates $\tilde{\rho}$ defined by $ \tilde{\rho}(\gamma)= (x^{\epsilon(\gamma)}, \rho(\mu^{-\epsilon(\gamma)}\gamma)).$}
Given $\chi: H \to \mathbb{Z}_k$,  define $\Phi: \pi_1(X) \xrightarrow{\tilde{\rho} }\mathbb{Z} \ltimes H \to GL_2(\mathbb{Q}(\xi_k)[t^{\pm1}])$ as the composition of $\tilde{\rho}$ with the map
$
(x^j, v) \mapsto
 \left[
 \begin{array}{cc}
 0& 1 \\
 t & 0 
 \end{array}
 \right]^j
 \left[
 \begin{array}{cc}
 \xi_k^{\chi(v)}&0\\
 0 & \xi_k^{-\chi( v)}
 \end{array}
 \right]
$.
Then we have the following.
\begin{thm}[\cite{HKL10}]\label{covercomp}
Let $X, X_2, \epsilon, \chi, \rho,$ and $\Phi$ be as above, where
\begin{itemize}
\item
  $\mathbb{Q}(\xi_k)[t^{\pm1}]$ has a $(\mathbb{Q}(\xi_k)[t^{\pm1}], \mathbb{Z}[\pi_1(X_2)])$-bimodule structure with right action defined by $p(t) \cdot \gamma= \xi_k^{\chi\cdot \rho(\gamma)} t^{\epsilon_2(\gamma)} p(t)$. 
\item $(\mathbb{Q}(\xi_k)[t^{\pm1}])^2$ has a $(\mathbb{Q}(\xi_k)[t^{\pm1}], \mathbb{Z}[\pi_1(X)])$-bimodule structure with right action defined by $\Phi: \pi_1(X) \to GL_2(\mathbb{Q}(\xi_k)[t^{\pm1}])$. 
\end{itemize}
The corresponding twisted homology groups $H_1(X_2,\mathbb{Q}(\xi_k)[t^{\pm1}])$ and $H_1(X, (\mathbb{Q}(\xi_k)[t^{\pm1}])^2)$ are isomorphic as $\mathbb{Q}(\xi_k)[t^{\pm1}]$-modules, and so the corresponding twisted Alexander polynomials are equal as well. 
\end{thm}

In practice, we define $\rho$ implicitly by constructing a map $\tilde{\rho}: \pi_1(X) \to \mathbb{Z} \ltimes H$ sending a Wirtinger generator $x_i$ to $(x, v_i)$ such that our preferred meridian $\mu$ is sent to $(x,0)$.  The Wirtinger relation $x_j x_i x_j^{-1}=x_k$  implies that we must have $(1-t) \cdot v_j+ t\cdot v_i= v_k$ in $H= \mathbb{F}_k[t]/\langle t+1 \rangle.$ However, since $t+1=0$  this relation reduces to $v_i+ v_k= 2 v_j$. 
We also need a choice of $\chi: H \to \mathbb{Z}_k$; since $H$ is one-dimensional over $\mathbb{Z}_k$, all nontrivial choices are essentially the same and so we take $\chi(1)=1$.

Finally, we need Wada's computationally powerful group-theoretic description of twisted Alexander polynomials, translated to the current context by Herald, Kirk, and Livingston \cite{Wada94, HKL10}. 
Suppose that $\pi=\pi_1(X)= \langle x_1, \dots, x_{s+1} : r_1, \dots, r_{\frak{s}} \rangle$, where $X=X(K)$ is homotopy equivalent to a CW complex with a single 0-cell, $(s+1)$ 1-cells, and $\frak{s}$ 2-cells. Let $\frac{\partial r_i}{\partial x_j}$ denote the Fox derivative of $r_i$ with respect to $x_j$. 
Let $\rho: \pi \to GL_n(\mathbb{F})$ and $\epsilon: \pi \to \mathbb{Z}= \langle t \rangle$ be nontrivial. Define  $F$ to be the composition $F: \mathbb{Z}[\langle x_1, \dots, x_{s+1} \rangle] \twoheadrightarrow \mathbb{Z}[\pi] \xrightarrow{\epsilon \otimes \rho} M_n(\mathbb{F}[t^{\pm1}]).$
Then the twisted chain complex $C_*=C_*(X, \mathbb{F}[t^{\pm1}]^n)$ has $C_2=(\mathbb{F}[t^{\pm1}]^n)^{\frak{s}}$, $C_1= (\mathbb{F}[t^{\pm1}]^n)^{s+1}$, and $\partial_2: C_2 \to C_1$ given by the block matrix $\left[  F(\frac{\partial r_i}{\partial x_j})  \right]_{\frak{s},s+1}$.

\begin{thm}[\cite{Wada94}, \cite{KL99a}] \label{foxder}
With the setup above, there is some $k$ such that $F(x_k-1)$ has nonzero determinant. 
Let $p_k: (\mathbb{F}[t^{\pm1}]^n)^{s+1} \to  (\mathbb{F}[t^{\pm1}]^n)^{s}$ be the projection with kernel the $k^{th}$ copy of $\mathbb{F}[t^{\pm1}]^n$.  
Define $Q_k \in \mathbb{F}[t^{\pm1}]$ to be the greatest common divisor of the 
$ns \times ns$ subdeterminants of the matrix for $p_k \circ \partial_2: (\mathbb{F}[t^{\pm1}]^n)^{\frak{s}}\to  (\mathbb{F}[t^{\pm1}]^n)^{s}$. 
Then, when $H_1(X, \mathbb{F}[t^{\pm1}]^n)$ is torsion,
\[\Delta(X, \mathbb{F}[t^{\pm1}]^n)= Q_k \frac{ \Delta_0(X)}{\det(F(x_k-1))}\]
\end{thm} 

In our case, we will have a generator $\mu=x_k$ in $\pi_1(X)$ with $\chi(x_k)=0$ and $\epsilon(x_k)=1$,  so $\Delta_0(X)=1$. In addition, we will choose $\tilde{\rho}$ so that for some generator $x_k$, we have $\det(F(x_k-1))=1-t$. Finally, we will work with a Wirtinger presentation, which has deficiency one (i.e., $\frak{s}=s$) and hence eliminates the need to take greatest common divisors.  So we will have $\Delta(X, \mathbb{F}[t^{\pm1}]^n)= \det F(Z) (1-t)^{-1}$, where $Z$ is obtained from $\left[ \frac{\partial r_i}{\partial x_j} \right]_{s,s+1} $ by deleting the block column corresponding to $x_k$.

\section{Results}
We have the following set-up. Let $K=K_{p,q}$ be a 2-bridge knot with Wirtinger presentation   $\pi_1(X)=\langle x_1, \dots, x_{s+1} | \, r_1, \dots, r_{s} \rangle$.  Suppose $p=m^2$ and let $k$ be a prime dividing $m$. 
Let $\tilde{\rho}: \langle x_1, \dots, x_{s+1} | \, r_1, \dots, r_{s} \rangle \to \mathbb{Z} \ltimes \mathbb{F}_k$ be any map such that
$\tilde{\rho}(x_i)= (x,v_i)$ for $i=1, \dots, s$ ,  $\tilde{\rho}(x_{s+1})=(x,0)$, and such that whenever $x_j x_i x_j^{-1} x_l^{-1}$ is a relation then we have that $2v_j= v_i + v_l$.\footnote{That is, $\tilde{\rho}$ is a homomorphism of the desired form.}
Let $\Phi: \pi_1(X) \to GL_2(\mathbb{Q}(\xi_k)[t^{\pm1}])$ be defined by
\[ x_i \mapsto (x, v_i) \mapsto 
\left[
\begin{array}{cc}
0&1\\
t&0\\
\end{array}
\right]
\left[
\begin{array}{cc}
\xi_k^{v_i}&0\\
0&\xi_k^{-v_i}\\
\end{array}
\right]
= 
\left[
\begin{array}{cc}
0&\xi_k^{-v_i}\\
t\xi_k^{v_i}&0\\
\end{array}
\right]
,\] and let $F_{\Phi}$ be the natural extension $\mathbb{Z}[\pi_1(X)] \to M_2(\mathbb{Q}(\xi_k)[t^{\pm1}])$. 
If $K$ is topologically slice, then
\[ \widetilde{\Delta^\Phi_K}(t)=(t-1)^{-2} \det F_{\Phi} \left(\left[ \frac{\partial r_i}{\partial x_j} \right]_{s,s}\right) \in \mathbb{Q}(\xi_k)[t^{\pm1}]\]
must factor as a norm in $\mathbb{Q}(\xi_k)[t^{\pm1}]$.

Note that the computation of $\widetilde{\Delta^\Phi_K}(t)$ as described above is   easy to implement on a computer. To obstruct  the topological sliceness of $K_{p,q}$ we can assume, switching $(p,q)$ with $(p,p-q)$ if necessary, that $q$ is even and so $p/q$ has an even continued fraction expansion. There is a straightforward formula for the Wirtinger presentation of $\pi_1(X(K_{p,q}))$ in terms of this even continued fraction expansion, and we obtain $\tilde{\rho}$ by solving a simple system of linear equations over $\mathbb{F}_k$. The twisted Alexander polynomial is then obtained via a simple computation; the only non-algorithmic part comes in showing that a particular $\widetilde{\Delta^\Phi_K}(t)$ does not factor as a norm in $\mathbb{Q}(\xi_k)[t^{\pm1}]$.

\begin{exl}\label{225}
When $K=K_{225,94}$ we have continued fraction expansion $[2,2,2,-6,-2,2]$ and Alexander polynomial $(3t^3-6t^2+5t-1)(t^3-5t^2+6t-3)$. It is straightforward to check that $K$ is algebraically slice and that all prime-power Casson-Gordon signature invariants vanish. 
However, there are Casson-Gordon signatures that obstruct $K$ from being ribbon, and Lisca's results show that $K$ is not even smoothly slice. We can show that $K$ is not topologically slice via the computation of a single twisted Alexander polynomial, corresponding to $k=5.$ (It is perhaps interesting to note that the twisted Alexander polynomial corresponding to $k=3$ factors as a norm even in $\mathbb{Q}[t^{\pm1}]$.)

The reduced twisted Alexander polynomial corresponding to $k=5$ is given by $\widetilde{\Delta^\Phi_K}(t)=(2+\xi_5^2 + \xi_5^3) (t^4+1) - (18+11(\xi_5^2+\xi_5^3))(t^3+t) + (34+21(\xi_5^2+\xi_5^3)) t^2$.
Note that since $\xi_5^2 + \xi_5^3= \frac{1}{2} (-1 - \sqrt{5})$, we have that, up to multiplication by units,
\[ 
\widetilde{\Delta^\Phi_K}(t)=
(3- \sqrt{5}) (t^4+1)- (25-11\sqrt5)(t^3+1) + (47-21 \sqrt5)t^2. 
\]
To show that $K_{225,94}$ is not slice, we must obstruct this polynomial from factoring  as a norm in $\mathbb{Q}(\xi_5)[t^{\pm1}]$. 
Consider the Galois conjugate $g(t) = (3+ \sqrt{5}) (t^4+1)- (25+11\sqrt5)(t^3+1) + (47+21 \sqrt5)t^2$. Note that any factorization of $\widetilde{\Delta^\Phi_K}(t)$ in $\mathbb{Q}(\xi_5)[t^{\pm1}]$ induces a corresponding factorization of $g(t)$, so it suffices to show that $g(t)$ is not a norm over $\mathbb{Q}(\xi_5)$. In fact, $g(t)$ has four distinct real roots and so it is enough to obstruct $g(t)$ from factoring as a norm over $\mathbb{Q}(\xi_5) \cap \mathbb{R}= \mathbb{Q}(\sqrt{5})$. 
So suppose that there are $\lambda, a, b, c \in \mathbb{Q}(\sqrt{5})$ such that 
$
g(t) 
= \lambda (at^2+bt+c)(ct^2+bt+a)$;
that is, such that
$\lambda a c = 3+ \sqrt{5}$,  $\lambda(a+c)b = - 25- 11 \sqrt{5}$, and $ \lambda(a^2+ b^2 + c^2)= 47+21\sqrt{5}$. 
This reduces to solving 
\[\frac{(a+c)b}{ac}=-5-2 \sqrt{5} \text{ and } \frac{a^2+b^2+c^2}{ac}=9+4\sqrt{5}  \text{ for } a, b, c \in \mathbb{Q}(\sqrt{5}). 
\]
It is straightforward to check using a computer algebra system that this has no solutions.
\end{exl}

\begin{exl} \label{cgfake}
We say $K_{m^2, q}$ is CG- fake slice if all prime-power Casson-Gordon signature invariants vanish but $K$ is not ribbon (or, equivalently by \cite{Lisca07}, not smoothly slice). 
The following table gives a count, for each $m$, of how many $K_{m^2,q}$ are CG-fake slice (counting $K$ and $-K$ as a single entry). We omit $m$ which are prime powers, since our computations agree with the conjecture that in this case CG signatures exactly detect smooth sliceness. These computations were done in Sage.

\begin{table}[htdp]
\begin{center}
\begin{tabular}{|c|c|c|}
m & Number of CG-fake slice $K_{m^2,q}$  & Number with $\Delta_K(t)$ a norm \\
\hline
$3 \cdot 5$ & 2&1 \\
$3 \cdot 7$& 3&0 \\
$3 \cdot 11$&3&0\\
$5 \cdot 7$ &10&2 \\
$3 \cdot 13$ & 5 & 0 \\
$3^2 \cdot 5$ &3 & 0 \\
$3\cdot 17$ & 5 & 0 \\
$5 \cdot 11$ & 16 & 2 \\
$3 \cdot 19$ & 3 & 0
\end{tabular}
\end{center}
\label{default}
\caption{Failure of Casson-Gordon signatures and Alexander polynomials to obstruct smooth sliceness}
\end{table}%
\end{exl}

\begin{exl}\label{1225}
The next knot we are led to consider is $K=K_{1225, 466}$.  $K$ has even continued fraction expansion $[2,2,-2,-2,-4,4,2,-2]$ and Alexander polynomial $(t^4-6t^3+13t^2-11t+4)(4t^4-11t^3+13t^2-6t+1)$. Again, $K$ is algebraically slice and has all prime-power CG signature invariants trivial, but is not smoothly slice. 
The twisted Alexander polynomial corresponding to $k=7$ is
\begin{align*}
\widetilde{\Delta^\Phi_K}(t)=& (8+4(\xi^3+\xi^4))(t^6+1)-(81+48(\xi^3+\xi^4)-16(\xi^2+\xi^5))(t^5+t)\\
&+ (287+189(\xi^3+\xi^4)-45(\xi^2+\xi^5)+27(\xi+\xi^6))(t^4+t^2)\\
&- (300+160(\xi^3+\xi^4)-188(\xi^2+\xi^5)-75(\xi+\xi^6))t^3.\\
\end{align*}
To show that this polynomial does not factor as a norm in $\mathbb{Q}[\xi_7]$, we use the following extension of Gauss' Lemma from Herald, Kirk, and Livingston. 
\begin{lem} \label{factoring} \cite{HKL10}
Let $k$ and $r$ be primes such that $r=nk+1$ for some $n \in \mathbb{N}$. Let $b \in \mathbb{Z}_r$ be a nontrivial $k^{th}$ root of 1, and let $\phi: \mathbb{Z}[\xi_k] \to \mathbb{Z}_r$ be the ring homomorphism sending $1$ to $1$ and $\xi_k$ to $b$. 
Let $p(t) \in \mathbb{Z}[\xi_k](t)$ be a degree $2m$ polynomial, such that $\phi(p(t))$ also has degree $2m$. 

If $p(t)$ is a norm in $\mathbb{Q}[\xi_k](t)$, then $\phi(p(t))$ factors as the product of two degree $m$ polynomials in $\mathbb{Z}_r[t]$. 
\end{lem}

In this case, we take $k=7$, $r=29=4*7+1$, and $b=16 \in \mathbb{Z}_{29}$. Let $\phi: \mathbb{Z}[\xi_7] \to \mathbb{Z}_{29}$ be defined as above with $1 \mapsto 1$ and $\xi_7 \mapsto 16$. Then $\phi\left(\widetilde{\Delta^\Phi_K}(t)\right)= 20(1+6t+t^2)(1+16t+6t^2+16t^3+t^4)$ is still degree 6 and has a $\mathbb{Z}_{29}$-irreducible degree 4 factor. So, by Lemma \ref{factoring}, $\widetilde{\Delta^\Phi_K}(t)$ is  not a norm over $\mathbb{Q}[\xi_7]$ and hence $K$ is not topologically slice. 
\end{exl}

Note that the above arguments obstructing $\widetilde{\Delta^\Phi_K}(t)$ from factoring as a norm in the appropriate field are quite ad hoc, and there is no reason to believe that either would necessarily be effective for a larger class of 2-bridge knots. In fact, each argument fails to work for the other example. This is emphasized even more by our computations for $K_{1225, 496}$. The reduced twisted Alexander polynomial for $K$ corresponding to a nontrivial character to $\mathbb{Z}_5$ factors as a norm.  While the polynomial corresponding to a nontrivial character to $\mathbb{Z}_7$ is not obviously a norm, both of the strategies used in Examples \ref{225} and \ref{1225} fail to obstruct such a factorization. 

\section*{Acknowledgments}
I would like to thank my advisor Cameron Gordon for suggesting this problem to me, as well as for his encouragement and advice.  


\bibliography{rational}

\begin{thebibliography}{HKL10}

\bibitem[CG86]{CG86}
Andrew Casson and Cameron Gordon.
\newblock Cobordism of classical knots.
\newblock In Guillou and Marin, editors, {\em A la recherche de la {Topologie}
  perdue}, volume~62 of {\em Progress in Mathematics}. Birkh{\"a}user Boston,
  1986.

\bibitem[EL09]{EL09}
Michael Eisermann and Christoph Lamm.
\newblock For which triangles is {P}ick's formula almost correct?
\newblock {\em Experimental Mathematics}, 18:187--191, 2009.

\bibitem[HKL10]{HKL10}
Chris Herald, Paul Kirk, and Charles Livingston.
\newblock Metabelian representations, twisted {A}lexander polynomials, knot
  slicing, and mutation.
\newblock {\em Math Z.}, 265(4):925--949, 2010.

\bibitem[KL99]{KL99a}
Paul Kirk and Charles Livingston.
\newblock Twisted {A}lexander invariants, {R}eidemeister torsion, and
  {C}asson-{G}ordon invariants.
\newblock {\em Topology}, 38(3):635--661, 1999.

\bibitem[Lis07]{Lisca07}
Paolo Lisca.
\newblock Lens spaces, rational balls and the ribbon conjecture.
\newblock {\em Geom.Topol.}, 11:429--473, 2007.

\bibitem[Wad94]{Wada94}
Masaaki Wada.
\newblock Twisted {A}lexander polynomial for finitely presentable groups.
\newblock {\em Topology}, 33(2):241--256, 1994.

\end{thebibliography}

\end{document}